\documentclass[12pt,psamsfonts]{article} 

\usepackage{a4,fullpage,amssymb,epsf,psfrag,times}
\usepackage{graphicx}
\usepackage{amsmath}
\usepackage{amsfonts}
\usepackage{enumerate}
\usepackage{latexsym}
\usepackage{epsfig}
\usepackage{amsthm}  
\usepackage{amsmath}
\usepackage{amssymb} 
\usepackage{latexsym}
\usepackage{epsfig} 
\usepackage{amssymb,latexsym}
\usepackage[all]{xy}
\usepackage{graphicx}
\usepackage{amsmath}
\usepackage{amsfonts}
\usepackage{enumerate}
\usepackage{latexsym}
\usepackage{epsfig}
\usepackage{amsthm}  
\usepackage{amsmath}
\usepackage{amssymb} 
\usepackage{latexsym}
\usepackage{epsfig}

\makeatletter\@addtoreset{equation}{section}\makeatother
\makeatletter\@addtoreset{table}{section}\makeatother

\newtheorem{theorem}{Theorem}[section]
\newtheorem{prop}[theorem]{Proposition}
\newtheorem{lemma}[theorem]{Lemma}

\newtheorem{cor}[theorem]{Corollary}

\newenvironment{remark}{\refstepcounter{theorem}\par\medskip\noindent{\bf Remark~\thetheorem~~}}{\unskip\nobreak\hfill\hbox{ $\oslash$}\par\bigskip}

\newenvironment{question}{\refstepcounter{theorem}\par\medskip\noindent{\bf Question~\thetheorem~~}}{\unskip\nobreak\hfill\hbox{ $\oslash$}\par\bigskip}

\newenvironment{example}{\refstepcounter{theorem}\par\medskip\noindent{\bf Example~\thetheorem~~}}{\unskip\nobreak\hfill\hbox{ $\oslash$}\par\bigskip}

\newenvironment{definition}{\refstepcounter{theorem}\par\medskip\noindent{\bf Definition~\thetheorem~~}}{\unskip\nobreak\hfill\hbox{ $\oslash$}\par\bigskip}

\begin{document}

\title{Toric symplectic ball packing}
\date{}
\author{Alvaro Pelayo}

\maketitle

\begin{abstract}
We define and solve the toric version of the symplectic ball packing problem, 
in the sense of listing all $2n$\--dimensional symplectic--toric
manifolds which admit a perfect packing by balls embedded
in a symplectic and torus equivariant fashion.

In order to do this we first describe a problem in geometric--combinatorics
which is equivalent to the toric symplectic ball packing problem.
Then we solve this problem using arguments from Convex Geometry
and Delzant theory.

Applications to symplectic blowing--up are also presented,
and some further questions are raised in the last section.
\end{abstract}

\section{The Main Theorem}

Loosely speaking, the ``symplectic packing problem''
asks how much of the volume of a symplectic manifold $(M,\, \sigma)$
may be filled up with disjoint embedded \emph{open} symplectic balls.
A lot of progress on this and intimately related questions
has been made by a number of authors, among them Biran \cite{B0}, \cite{B1/2}, \cite{B3}, 
McDuff--Polterovich \cite{MP}, Traynor \cite{T} and Xu \cite{Xu}.

Several authors have made progress on directly related questions, like the
topology of the space of symplectic ball embeddings, among them
McDuff \cite{M1}, Biran \cite{B1} and most recently Lalonde--Pinsonnault \cite{LP} (the equivariant version of this question
was studied in \cite{P}). 

Despite the fact that these significant contributions have appeared in recent years, 
the symplectic packing problem remains largely not understood; 
for more details and a survey of known results see 
the paper by Biran \cite{B2} and for some nice examples see \cite{MMT}. 
Outstanding progress has been made in dimension four (see for example \cite{B0}, \cite{MP}, \cite{T}), 
but nothing is understood in dimension six or above. 
The underlying reason for this dimensional barrier is that 
the techniques used by the previous
authors are unique to dimension four and do not extend to higher dimensions.
For more details see Section 1 of \cite{P} or Section 3 of Biran's paper \cite{B3}.

The present paper is devoted to the study of a particular case of the symplectic packing problem, 
the \emph{torus equivariant case}. In this case both the symplectic manifold $M$ and the standard
open symplectic ball $(\mathbb{B}_r,\, \sigma_0)$ in $\mathbb{C}^n$ are equipped with a Hamiltonian action of
an $n$\--dimensional torus $\mathbb{T}^n$, and the symplectic embeddings 
of this ball into $M$ that we consider are equivariant with respect to these
actions. 

Our main result is Theorem \ref{mt}, which provides the list
of symplectic--toric manifolds which admit a full packing by balls embedded
in such a way, i.e. we prove existence and a uniqueness of such manifolds. Our proofs rely on the discovery by Delzant \cite{D} that symplectic--toric
manifolds are classified by their convex images under the momentum map. This allows us to solve an a priori symplectic--geometric
problem using techniques from Convex Geometry and Delzant theory (for a treatment of this 
theory see for example the book by Guillemin \cite{G3}).

\begin{figure} \label{AFF}
\begin{center}
\epsfig{file=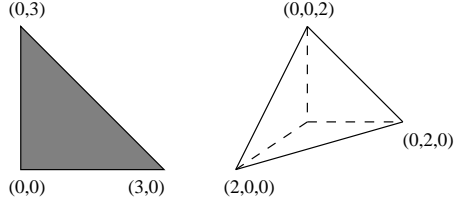}
\caption{Momentum polytope of $(\mathbb{CP}^2, \, 3 \cdot \sigma_{\textup{FS}})$ (left) and $(\mathbb{CP}^3, \, 2 \cdot \sigma_{\textup{FS}})$ (right).}
\end{center}
\end{figure}

\begin{definition} 
\label{Definition A}(\cite{C}).
A compact connected symplectic manifold $M=(M, \, \sigma)$ is a \emph{symplectic--toric manifold},
or a \emph{Delzant manifold}, if it is equipped with an effective and Hamiltonian
action of a torus of dimension half of the dimension of the manifold.
\end{definition}

Delzant manifolds come equipped with a momentum map $\mu^M \colon M \to \mathbb{R}^n$
which satisfies $\textup{i}_{\xi_M} \sigma=\textup{d} \langle \mu^M,\, \xi \rangle$
for all $\xi$ in the Lie algebra of the torus (see \cite{G3}), and where
$\xi_M$ is the vector field on $M$ induced by $\xi$ (the existence of $\mu$ 
may be taken as definition of Hamiltonian action).
$\mu^M$ carries $M$ to a convex polytope in $\mathbb{R}^n$ (here we
identify $\mathbb{R}^n$ with the dual of the Lie algebra of the torus; see Section 2 for details
on how we make this non--canonical identification),
and this polytope, which is called the \emph{momentum polytope of $M$}, determines $M$ up to equivariant
symplectomorphism (see Theorem \ref{Theorem I}). 

\begin{example} \label{Example B} (Projective Spaces). The projective space $(\mathbb{CP}^n,\, \lambda \cdot \sigma_{\textup{FS}})$  
equipped with a $\lambda$ multiple of  the Fubini--Study form
$
\sigma_{\textup{FS}}=\frac{1}{2 (\sum_{i=0}^n \bar{z}_i z_i)} \sum_{k=0}^n \sum_{j \neq k} (\bar z_j z_j \, \textup{d}z_k \wedge \textup{d} \bar{z}_k-
\bar z_j z_k  \, \textup{d}z_j \wedge \textup{d} \bar{z}_k)
$
and the rotational action of $\mathbb{T}^n$, 
$
(\textup{e}^{\textup{i} \theta_1},\ldots,\textup{e}^{\textup{i} \theta_n}) \cdot [z_0: \ldots:z_n]=[z_0:\, \textup{e}^{-2 \pi \textup{i} \theta_1}\, z_1: \ldots:\textup{e}^{-2 \pi \textup{i} \theta_n} \, z_n]
$,
is a $2n$--dimensional Delzant manifold with momentum map components
$\mu^{\mathbb{CP}^n}_k(z)=\frac{\lambda |z_k|}{\sum_{i=0}^n|z_i|^2}$
and whose momentum polytope equals the convex hull in $\mathbb{R}^n$ of $0$ and 
the scaled canonical vectors $\lambda e_1,\ldots,\lambda e_n$, see Figure 2.
\end{example}

\begin{example} \label{Example C}
The product $(\prod \mathbb{CP}^{n_i}, \, \prod \lambda_i \cdot \sigma_{\textup{FS}})$ is
a Delzant manifold of dimension $2 \, \sum n_i$.
\end{example}

\begin{example} \label {Example D} (Open Balls). The open symplectic ball $(\mathbb{B}_r,\, \sigma_0)$
in $\mathbb{C}^n$ with the $\mathbb{T}^n$ action by rotations
(component by component) has momentum map components
$\mu_k^{\mathbb{B}_r}(z)=|z_k|^2$
and momentum polytope
equal to the convex hull in $\mathbb{R}^n$ of $0$
and the scaled canonical basis vectors $r^2e_1,\ldots, 
r^2e_n$.

The momentum polytope of $(\mathbb{CP}^n,\, \lambda \cdot \sigma_{\textup{FS}})$
equals the momentum polytope of $\mathbb{B}_{\sqrt{\lambda}}$ \emph{minus} the face 
of the momentum polytope of the former which faces the origin, which does
not belong to the momentum polytope $\Delta^{\mathbb{B}_r}$ of $\mathbb{B}_r$.
$\Delta^{\mathbb{B}_r}$ is an integral simplex in the sense of Definition \ref{setup}.
\end{example}

\begin{definition} \label{Definition F} An embedding $f$ of the $2n$\--ball $\mathbb{B}_{r}$ into a 
$2n$--dimensional Delzant manifold $M$ is \emph{equivariant} if there exists
an automorphism $\Lambda$ of $\mathbb{T}^n$ such that the  diagram 
\begin{eqnarray} 
\xymatrix{ \ar @{} [dr] |{\circlearrowleft}
\mathbb{T}^{n} \times \mathbb{B}_r  \ar[r]^{\Lambda \times f}      \ar[d]^{ \mathbf{\cdot} }  &  \mathbb{T}^{n} \times M 
                  \ar[d]^{\psi}   \\
                   \mathbb{B}_r  \ar[r]^f   &       M    } \nonumber
\end{eqnarray} 
commutes, where $\psi$ is a fixed effective and Hamiltonian $\mathbb{T}^n$\--action on $M$
and $\mathbf{\cdot}$ denotes the standard action by rotations
on $\mathbb{B}_r$ \textup{(}component by component\textup{)}. In this case we say
that $f$ is a $\Lambda$--\emph{equivariant embedding}.
\end{definition}

Next we give a precise notion
of ``perfect equivariant symplectic ball packing''. In the following definition we use the term ``family of maps''
to mean a ``collection of maps'' in a set-theoretical fashion, i.e. we do not assume that this collection of
maps has any additional structure.

\begin{figure} \label{figure4}
\begin{center}
\epsfig{file=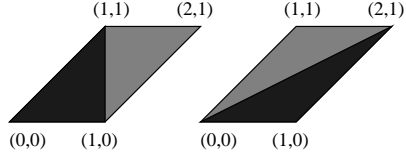}
\caption{A manifold equivariantly symplectomorphic to $(\mathbb{CP}^1 \times \mathbb{CP}^1, \, \sigma_{\textup{FS}} \oplus \sigma_{\textup{FS}})$
packed by two equivariant symplectic balls of radius $1$
in two different ways (see Lemma \ref{newlemma} for an explanation).}
\end{center}
\end{figure}

\begin{definition}\label{4.1}
We define the real--valued mapping $\Omega$ on the space of $2n$--dimensional Delzant manifolds by
$
\Omega(M):=(\textup{vol}_{\sigma}(M))^{-1} 
\, \, \sup_{\mathcal{E} \in \mathfrak{F}} \sum_{f \in \mathcal{E}} \textup{vol}_{\sigma}(f(\mathbb{B}_{r_{f}}))$,
where $\mathfrak{F}$ is the set of families $\mathcal{E}$ of equivariant symplectic ball embeddings such that
if $f,\, g \in \mathcal{E}$ then $f(\mathbb{B}_{r_f}) \cap g(\mathbb{B}_{r_g})=\emptyset$, 
and $r_{f} \ge 0$ for all $f \in \mathcal{E}$.
We say that $M$ \emph{admits a perfect equivariant and symplectic ball packing} if 
there exists a family $\mathcal{E}_0 \in \mathfrak{F}$ such that
$\Omega(M)=1$ at $\mathcal{E}_0$.
\end{definition}

A version of Definition \ref{4.1} in the ``general'' symplectic case, as well as the general symplectic
packing problem, were introduced by McDuff and Polterovich in \cite{MP}. They denote
$\Omega(M)$ by $\textup{v}(M,k)$ where $k$ is the (a priori fixed) number of balls that
we are embedding in $M$. 
McDuff and Polterovich consider that all balls have the same (a priori
fixed) radius $r>0$. This is in contrast to Definition \ref{4.1} above, where neither
the number of balls nor the radius are fixed. 

\begin{theorem} \label{mt}
Let $M$ be a $2n$--dimensional Delzant manifold. Then $M$ admits a perfect equivariant symplectic ball packing if and only if there exists $\lambda>0$ such that
\begin{enumerate}
\item
($n=2$) $M$ is equivariantly symplectomorphic to either $(\mathbb{CP}^{2}, \, \lambda \cdot \sigma_{\textup{FS}})$ or 
the product $(\mathbb{CP}^1  \times \mathbb{CP}^1, \, \lambda \cdot(\sigma_{\textup{FS}} \oplus \sigma_{\textup{FS}}))$.
\item
($n \neq 2$) $M$ is equivariantly symplectomorphic to $(\mathbb{CP}^{n},\, \lambda
\cdot \sigma_{\textup{FS}})$. 
\end{enumerate}
Equivalently,  $(\mathbb{CP}^{n},\, \lambda \cdot \sigma_{\textup{FS}})$ and 
$(\mathbb{CP}^1  \times \mathbb{CP}^1, \, \lambda \cdot(\sigma_{\textup{FS}} \oplus \sigma_{\textup{FS}}))$
with  $n \ge 1$ and $\lambda>0$ are the only Delzant manifolds manifolds which admit a perfect equivariant and symplectic ball packing. 
\end{theorem}

The proof of Theorem \ref{mt} 
in any dimension follows from the abstract combinatorial structure of Delzant polytopes, the abstract notion of convexity, and its properties. In addition, a number of figures are presented along with the proof to suggest some intuition of the solution.

In Section 3 we analyze in how many different ways the spaces which appear in Theorem \ref{mt} may be packed -- 
this is summarized in Proposition \ref{tt} below; the \emph{existence statement} in  Proposition \ref{tt}  is a simple construction, c.f. Remark \ref{Remark G}, Figure 2, Figure 3,
and it is independent from the proof of the \emph{uniqueness statement} in Theorem \ref{mt}, this last one being the part which occupies most of Section 3. 

The \emph{existence part} of the statement of Theorem \ref{mt}  follows from Proposition \ref{tt}. If Proposition \ref{tt}
is assumed, Theorem \ref{mt} becomes a uniqueness theorem.

\begin{prop} \label{tt}
For all $\lambda>0$ and $n \ge 1$, the complex projective space $(\mathbb{CP}^n,\, \lambda \cdot \sigma_{\textup{FS}})$ 
may be perfectly packed by one equivariant symplectic ball, and it
may not be perfectly packed by two or more equivariant symplectic balls for $n\geq 2$.

If $n=1$, $\lambda>0$, $(\mathbb{CP}^{n},\, \lambda \cdot \sigma_{\textup{FS}})$ 
may be perfectly packed only by one 
or two equivariant symplectic balls, and there is a one parameter family of packings by two equivariant symplectic
balls, c.f. Figure 3. For 
all $\lambda>0$, the $4$--dimensional Delzant manifold $(\mathbb{CP}^1 \times
\mathbb{CP}^1, \, \lambda \cdot (\sigma_{\textup{FS}} \oplus \sigma_{\textup{FS}}))$ 
may only be perfectly packed by two equivariant symplectic balls, 
and this in precisely two distinct ways, c.f. Figure 2. 
\end{prop}

\begin{figure}
\begin{center}
\epsfig{file=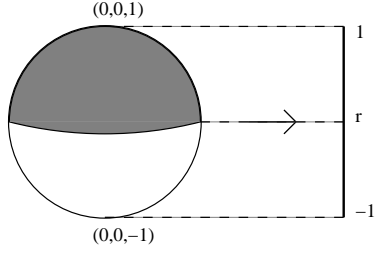}   
\label{AF6}
\caption{A one parameter family of perfect packings of $(S^2, \, \frac{1}{2} \cdot \textup{d}\theta \wedge \textup{d} h) \simeq (\mathbb{CP}^1,\, 2 \cdot \sigma_{\textup{FS}})$ by 
a ball of radius $0 \le R \le \sqrt{2}$ and a ball of radius $\sqrt{2-R^2}$, see Remark \ref{weird}. In the figure $\textup{r}=0$ and $R=1$. In general, $\textup{r}=1-R^2$.}
\end{center}
\end{figure}

\begin{remark} \label{weird}
\normalfont
We identify the $2$\--sphere of radius $R$ equipped with the standard area form $\textup{d} \theta \wedge \textup{d}h$ with $(\mathbb{CP}^1, \, 4R^2 \cdot \sigma_{\textup{FS}})$, where $\sigma_{\textup{FS}}$ is the Fubini--Study form on $\mathbb{CP}^1$. Under the conventions
that we use throughout the paper, the momentum map for 
$(S^2, \, \frac{1}{2} \cdot \textup{d}\theta \wedge \textup{d}h)$ is
equal to $(\theta, \, h)  \mapsto h$, and the momentum polytope is $[-1,\,1]$. In the literature it seems (far) more common to have $[-1,\,1]$ as momentum polytope for $S^2$ with the standard area form
$\textup{d}\theta \wedge \textup{d}h$; we do not follow this convention
in order to have a simpler expression for $\mu^M(f(\mathbb{B}_r))$ in Remark \ref{Remark G}.
Notice that the area of $(S^2, \, \frac{1}{2} \cdot \textup{d}\theta \wedge \textup{d}h)$ is $2 \, \pi$, while
the length of the associated momentum polytope is $2$.
\end{remark}

The paper is divided into five sections: in Section 2 we describe
a problem in geometric--combinatorics equivalent to the
toric symplectic ball packing problem; in Section 3
we solve it; in Section 4 we relate our results to the theory
of blowing up; we end by raising some further questions
in Section 5.

\section{From Symplectic Geometry to Combinatorics}
Let $M$ be a Delzant
manifold of dimension $2n$, and denote by $\mathbb{B}_{r}$ the $2n$--dimensional
ball in $\mathbb{C}^n$ equipped with the restriction of the standard symplectic form $\sigma_{0}$,
and with the standard action by rotations of the $n$--torus $\mathbb{T}^{n}$ (see Definition \ref{Definition A}; for the main properties of Delzant
manifolds see Section 1 of  \cite{P}, or for more details \cite{C}, \cite{D}, \cite{G3}).

Recall that the main feature that makes the study of symplectic manifolds equipped with torus actions 
richer than the study of generic symplectic manifolds is the existence of the
momentum map
$
\mu^{M} \colon M \to \mbox{Lie}(\mathbb{T}^n)^*
$
whose image $\Delta^M$ is a convex polytope, called the \emph{momentum polytope
of $M$}, as shown independently by Atiyah and Guillemin--Sternberg \cite{A}, \cite{G1}. 
The momentum map is unique up to addition of a constant in $(\textup{Lie}(\mathbb{T}^n))^*$,
and it is in this sense that we say ``the'' momentum map instead of ``a'' momentum map.
Here we are identifying the Lie algebra $\mbox{Lie}(\mathbb{T}^n)$ and its dual 
$\mbox{Lie}(\mathbb{T}^n)^*$ with $\mathbb{R}^n$. We denote by $\chi(M)$ the Euler characteristic of $M$.

Since this identification is not canonical we need to specify 
the convention we adopt in this paper. This amounts to choosing
an epimorphism $ \mathbb{R} \to \mathbb{T}^1$ which we take to be
$x \mapsto e^{2\sqrt{-1}x}$. This epimorphism induces an isomorphism
between $\mbox{Lie}(\mathbb{T}^1)$ and $\mathbb{R}$ via
$\frac{\partial}{\partial x} \mapsto 1/2$ giving rise to a new isomorphism
$\mbox{Lie}(\mathbb{T}^n) \to \mathbb{R}^n$, $\frac{\partial}{\partial x_k} \mapsto 1/2\, e_k$, 
by canonically identifying $\mbox{Lie}(\mathbb{T}^n)$ with the product of $n$ copies
of $\mbox{Lie}(\mathbb{T}^1)$ (see Section 32 in \cite{Gnew} for more details).

\begin{remark}
\label{Remark G}
\normalfont
In Section 2 of \cite{P} we proved that if $f:\mathbb{B}_r \to M$
is a $\Lambda$--equivariant and symplectic embedding with $f(0)=p$
and $\mu^M(p)=x$, then the following diagram is commutative:
\begin{eqnarray} \label{ml} 
\xymatrix{ \ar @{} [dr] |{\circlearrowleft}
\Delta^{\mathbb{B}_r} \ar[r]^{(\Lambda^{t})^{-1}+x}      &  \Delta^M  \\
           \mathbb{B}_r  \ar[r]^f  \ar[u]^{\mu^{\mathbb{B}_r}}  &       M      \ar[u]_{\mu^M}}
\end{eqnarray}
It follows from diagram (\ref{ml}) that if $M$ is a $2n$--dimensional Delzant manifold and
$f$ is a symplectic $\Lambda$--equivariant embedding from $\mathbb{B}_r$
into $M$ with $f(0)=p$ and $\mu^M(p)=x$,
then the momentum image
$\mu^M(f(\mathbb{B}_r))$ equals the subset of $\mathbb{R}^{n}$ given by the convex hull of 
the points $x$ and $x+r^2 \, \alpha^{p}_{i}$, $i=1,\ldots,n$,
minus the convex hull of  $x+r^2 \, \alpha^{p}_{i}$, $i=1,\ldots,n$,
where the $\alpha_i^{p}$
are the characters of the isotropy representation of $\mathbb{T}^{n}$ on the tangent space at $p$ to $M$.
\end{remark}

In the case when $M$ is a Delzant manifold, the momentum polytope $\Delta^{M}$ of $M$
is the so called \emph{Delzant polyope of $M$}, and it satisfies specific properties
as we see from Definition \ref{Definition H} below, which is a purely combinatorial definition -- in 
this respect Delzant polytopes may be defined without reference to
Delzant manifolds.

\begin{figure}[htb]
\begin{center}
\epsfig{file=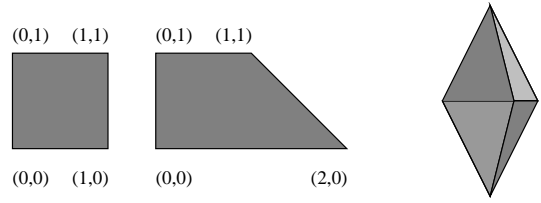}
\label{AF2}
\caption{Two Delzant polytopes (left) and a non-Delzant polytope (right).}
\end{center}
\end{figure}

\begin{definition}  [\cite{G3}, \cite{C}]  \label{Definition H}
A \emph{Delzant polytope $\Delta$ of dimension $n$ in $\mathbb{R}^{n}$}
is a simple, rational and smooth polytope. Here \emph{simple} means that there
are exactly $n$ edges meeting at each vertex of $\Delta$; \emph{rational} means that
the edges of $\Delta$ meeting at the vertex $x$ 
are rational in the sense that each edge is of the form
$x+t\, u_{i}$, $t \ge 0$, $u_i \in \mathbb{Z}^{n}$; \emph{smooth}
means that $u_1,\ldots,u_n$ may be chosen
to be a $\mathbb{Z}$\--basis of the integer lattice $\mathbb{Z}^{n}$.
\end{definition}

\begin{remark}
\normalfont
A similar class of polytopes, \emph{Newton polytopes},
have long been considered in algebraic geometry \cite{F}; the (only) difference 
between Newton and Delzant polyopes is that the former
are required to have all of their vertices lying
in the integer lattice $\mathbb{Z}^n$.

We learned from Yael Karshon that ``Delzant polytopes'' have also been refered to as
``non--singular'', ``torsion--free'' or ``unimodular'' by other authors.
\end{remark}

\begin{theorem}  [Delzant, \cite{D}]
\label{Theorem I}
Two Delzant manifolds are equivariantly symplectomorphic if and only if they
have the same Delzant polytope up to a transformation in $\textup{SL}(n,\,\mathbb{Z})$,
and translation by an element in $\mathbb{R}^n$.
For every Delzant polytope $\Delta$ there exists a Delzant manifold $M^{\Delta}$ whose 
momentum polytope is precisely $\Delta$.
\end{theorem}

\begin{definition}[Integral simplex] \label{setup}
If $\Upsilon$ is an $n$--dimensional simplex \textup{(}i.e. closed convex hull in $\mathbb{R}^n$
of $n+1$ linearly independent points\textup{)}, an \emph{open simplex with respect
to a vertex $x$ of} $\Upsilon$ is the convex region obtained
from $\Upsilon$ by removing the only face of $\Upsilon$
which does not contain $x$. A region $\Sigma \subset \mathbb{R}^n $ is
an \emph{open simplex} if there is a simplex $\Upsilon$ such that $\Sigma$
is an open simplex with respect to $x$ for some vertex $x$
of $\Upsilon$. If $\Sigma$ is an open simplex, then its closure in $\mathbb{R}^n$
is denoted by $\Sigma_{\textup{c}}$. We say that an $n$--dimensional open simplex $\Sigma \subset \mathbb{R}^n$
is \emph{integral} if:
\begin{enumerate}
\item
The  $\textup{SL}(n,\, {\mathbb{Z}})$\--length of each edge of $\Sigma$ is the same for all edges,
\item
$\Sigma$ has the same Euclidean volume as a simplex $d \, \Delta^{0}$ for some $d \ge 0$,
where $\Delta^{0}$ is the convex hull of $0$ and the canonical
basis vectors $e_{1},\ldots,e_{n}$.
\end{enumerate}
\end{definition}

\begin{remark}
\normalfont
If $\Sigma$ is a (closed) simplex of Euclidean volume equal to the Euclidean
volume of $d \, \Delta^{0}$, where $\Delta^{0}$ was defined in Definition \ref{setup},
then  all of the edges of $\Sigma$ meeting at a common vertex $x$ of $\Sigma$ 
have equal  $\textup{SL}(n\, \mathbb{Z})$\--length if and only if there exists a transformation
$A \in \textup{SL}(n\, \mathbb{Z})$ such that $A(d \, \Delta^{0})=\Sigma$. 
Similarly for open simplices -- notice that an open simplex has a unique vertex.
\end{remark}

\begin{figure}[htb]
\begin{center}
\epsfig{file=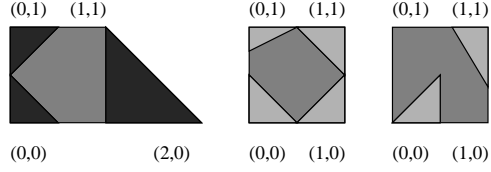}
\label{AF3}
\caption{The closures of a coherent (left) and two non-coherent families of simplices (right).}
\end{center}
\end{figure}

\begin{definition}[Coherent family of simplices] \label{cf}
Let $M$ a $2n$--dimensional Delzant manifold with Delzant
polytope $\Delta^{M}$. We say that a family $\mathcal{E}^{\Delta}$ of 
$n$--dimensional open simplices contained in $\Delta^{M}$ is \emph{coherent} if
for every $\Sigma \in \mathcal{E}^{\Delta}$ the following three properties are satisfied:
\begin{enumerate}
\item
The vertex $x_{\Sigma}$ of $\Sigma$ is a vertex of $\Delta^M$,
\item
Every $(n-1)$\--dimensional face
of $\Sigma$ is contained in $\partial (\Delta^M)$,
\item
$\Sigma$ is an integral simplex.
\end{enumerate}
The family $\mathcal{E}^{\Delta}_{\textup{c}}$ of closed simplices $\Sigma_{\textup{c}}$, where
$\Sigma \in \mathcal{E}^{\Delta}$, is called the \emph{closure of $\mathcal{E}^{\Delta}$}.
\end{definition}

\begin{remark} \label{2.9}
\normalfont
It follows from Definition \ref{cf} that coherent families  of \emph{disjoint} simplices 
contain at most $\chi(M)$ simplices, 
where $\chi(M)$ stands for Euler characteristic of $M$. Therefore the closure of
such a family contains at most $\chi(M)$ simplices.
\end{remark}

\begin{theorem} [Atiyah \cite{A}, Guillemin--Sternberg \cite{G1}]
\label{Theorem K} 
Let $M$ be a $2n$--dimensional symplectic manifold 
equipped with an effective and Hamiltonian action of
an $m$--dimensional torus $\mathbb{T}^m$. Then
the fibers $(\mu^M)^{-1}(x)$, where $x
\in \Delta^M$, of the momentum mapping $\mu^M$, are connected
subsets of $M$.
\end{theorem}

The following is Theorem 2.10 in Guillemin's book \cite{G3}, adapted to fit the 
conventions for $\mu^M \colon M \to \mathbb{R}^n$, introduced on the third paragraph of
the section.

\begin{theorem} [\cite{G3}] \label{guillemin}
The symplectic volume of a Delzant manifold $Q$ of dimension $2m$
is equal to $m! \, \pi^m$ times the Euclidean volume of its momentum polytope $\Delta^Q$.
\end{theorem}

Additionally, unlike in \cite{G3} we use
$\textup{vol}_{\sigma}(S)=\int_S \, \sigma^n$, i.e. we do not normalize the integral by dividing by $n!$.
Although Theorem \ref{guillemin} is stated only for Delzant manifolds, the result extends
to $\mathbb{B}_r \subset \mathbb{C}^n$.

\begin{cor} \label{st}
The symplectic volume of $\mathbb{B}_r$ is equal to $n! \, \pi^n$ times the Euclidean
volume of $\Delta^{\mathbb{B}_r}$. If $f \colon \mathbb{B}_r \to M$ is an equivariant
symplectic embedding, the symplectic volume of $f(\mathbb{B}_r)$ is equal to $n! \, \pi^n$ times the
Euclidean volume of $\mu^M(f(\mathbb{B}_r))$, and to $n! \, \pi^n$ times the
Euclidean volume of $\Delta^{\mathbb{B}_r}$.
\end{cor}

Both Theorem \ref{guillemin} and Corollary \ref{st} are particular versions of the
Duistermaat--Heckman theorem in \cite{Du}, or Section 2 in \cite{G3}.
We can now describe $\Omega$ in combinatorial terms; we denote
the Euclidean volume measure in $\mathbb{R}^n$ by $\textup{vol}_{\textup{euc}}$.

\begin{lemma} \label{newlemma}
Let $M$ a $2n$--dimensional Delzant manifold and let $\Omega$ be the mapping defined
in Definition \ref{4.1}. Let 
$\mathfrak{F}^{\Delta}$ be the set of coherent families $\mathcal{E}^{\Delta}$ of
pairwise disjoint simplices contained in $\Delta^{M}$. Then:
\begin{eqnarray} \label{mm} 
\Omega(M)=\frac{1}{\textup{vol}_{\textup{euc}} (\Delta^M)} \, \,  
\sup_{\mathcal{E}^{\Delta} \in \mathfrak{F}^{\Delta}} 
\sum_{\Sigma \in \mathcal{E}^{\Delta}} \textup{vol}_{\textup{euc}}(\Sigma).
\end{eqnarray}
Furthermore, $M$ admits a perfect equivariant and symplectic ball packing
if and only if there exist a coherent family $\mathcal{E}^{\Delta}$ of pairwise disjoint simplices contained in $\Delta^{M}$
and such that $\Delta^{M}=\bigcup_{\Sigma \in \mathcal{E}^{\Delta}}\Sigma_{\textup{c}}$.
\end{lemma}

\begin{proof}
Write $\nu$ for the right hand--side of (\ref{mm}).
It follows from Theorem \ref{Theorem K} that a pairwise disjoint family of equivariant symplectic ball
embeddings $\mathcal{E}_0=\{f_i\}$ gives rise to
a pairwise disjoint family $\{\Sigma_{i}:=\mu^{M}(f_{i}(\mathbb{B}_{r_{i}}))\} \subset \Delta^{M}$. 
By Remark \ref{Remark G}, each
$\Sigma_i$ is an integral simplex contained in $\Delta^M$ and the family $\mathcal{E}_0$ is coherent.
Without loss of generality we assume that the supremum in the formula given in Definition \ref{4.1} is achieved
at the family $\mathcal{E}_0$. Since $\mu^{M}(f_i(\mathbb{B}_{r_{i}}))=\Sigma_{i}$ and $f_{i}$ is symplectic and equivariant, by Corollary \ref{st}
$
\textup{vol}_{\sigma}(f_{i}(\mathbb{B}_{r_{i}}))=n! \, \pi^n \, \textup{vol}_{\textup{euc}}(\Sigma_{i}) 
$
and $\textup{vol}_{\sigma}(M)=n! \, \pi^n \, \textup{vol}_{\textup{euc}}(\Delta^M)$, which by plugging these
values into the formula of $\Omega$ in Definition \ref{4.1} evaluated at the family $\mathcal{E}_0$,
implies that
$\Omega(M) \le \nu$. One shows that $\Omega(M) \ge \nu$, by starting with a coherent
family of pairwise disjoint simplices and repeating this same argument. 

Now suppose that $M$ admits a perfect equivariant and symplectic ball packing.
Then by Definition \ref{4.1}, 
$
\textup{vol}_{\sigma}(M)=\sum_{i} \textup{vol}_{\sigma}(f_i(\mathbb{B}_{r_i}))
$
at certain family of embedded balls,
and in this case the equality
$
\Delta^{M}=\bigcup_{i} \overline{\mu^M(f_{i}(\mathbb{B}_{r_{i}}))}
$
holds, and by Remark \ref{Remark G} each momentum image $\mu^{M}(f_{i}(\mathbb{B}_{r_{i}}))$ is an integral simplex;
finally since the images $f(\mathbb{B}_{r_{i}})$ are pairwise disjoint,
by Theorem \ref{Theorem K} these simplices are pairwise disjoint. The converse is proved similarly.
\end{proof}

We will use Lemma \ref{newlemma} in Step 1 of the proof of Theorem \ref{mt}.

\section{Proof of Theorem \ref{mt} and Proposition \ref{tt}}

For clarity the proof is divided into six steps:

\vspace{2mm}

{\bf Step 1}. 
In this step we analyze the type of simplices which both form a coherent family as well as give
rise to a perfect equivariant and symplectic packing.

\begin{lemma} \label{firstlemma}
If the coherent family $\mathcal{E}^{\Delta}$ contains only one open simplex, then $M$ is equivariantly
symplectomorphic to $(\mathbb{CP}^n, \, \lambda \cdot \sigma_{\textup{FS}})$ for some $\lambda>0$. Otherwise
there exist at least two disjoint simplices in $\mathcal{E}^{\Delta}$, and every simplex in the family has exactly one face which is
not contained in $\partial (\Delta^M)$.
\end{lemma}

\begin{proof}
By Lemma  \ref{newlemma}, there exist a coherent family $\mathcal{E}^{\Delta}$ consisting of pairwise 
disjoint open simplices $\Sigma$ contained in $\Delta^{M}$ such that
$
\Delta^{M}=\bigcup_{\Sigma \in \mathcal{E}^{\Delta}} \Sigma_{\textup{c}}
$.
The coherence of $\mathcal{E}^{\Delta}$ implies that
for each $\Sigma \in \mathcal{E}^{\Delta}$ there exists a vertex 
$x \in \Delta^{M}$ such that there are exactly $n$ faces of $\Sigma_{\textup{c}}$ of dimension $n-1$ which 
contain $x$, and each of them is contained in the boundary $\partial(\Delta^M)$
of $\Delta^{M}$, which leaves only the $(n-1)$\--dimensional face of $\Sigma_{\textup{c}}$ which does not contain $x$, say the face
$\widehat F_{\Sigma}$, as possibly not contained in the boundary $\partial(\Delta^M)$ of $\Delta^{M}$.
Call $\mathcal{U}_{0}$ to the subfamily of $\mathcal{E}^{\Delta}$
consisting of those simplices such that one of their
$(n-1)$\--dimensional faces is not contained in
$\partial (\Delta^{M})$, and $\mathcal{U}_{1}$
to the subfamily of $\mathcal{E}^{\Delta}$ such that all simplices of
this subfamily have all of their $(n-1)$\--dimensional faces  contained in the boundary
$\partial (\Delta^{M})$.  Denoting by
$
\Theta_{i}=\bigcup_{\Sigma  \in \mathcal{U}_{i}} \Sigma_{\textup{c}} \textup{ ,   } i=0, \, 1,
$
the observation made in the previous paragraph implies
that $\Theta_{0} \cup \Theta_{1}=\Delta^{M}$ and $\Theta_0 \cap \Theta_1= \emptyset$.

Now we distinguish two cases, according to 
whether $\Theta_{1}=\emptyset$ or $\Theta_{1}\neq \emptyset$.
Let us first assume that $\Theta_{1} \neq \emptyset$; then there exists
an open simplex $\Sigma$ contained in $\Delta^{M}$ such that all of the $(n-1)$\--dimensional 
faces of $\Sigma_{\textup{c}}$ are contained in $\partial (\Delta^{M})$, and therefore $\partial (\Sigma_{\textup{c}})\subset \partial (\Delta^{M})$.
This being the case, it follows from the convexity of $\Delta^{M}$
that $\Sigma_{\textup{c}}=\Delta^{M}$. Since $\Sigma_{\textup{c}}=\Delta^M$, and by construction 
$\Sigma \in \mathcal{E}^{\Delta}$, where $\mathcal{E}^{\Delta}$ 
is a coherent family by Theorem \ref{Theorem I}, we conclude that $\Delta^{M}$ is the Delzant polytope of a Delzant manifold $M$
equivariantly symplectomorphic to the $n$\--dimensional 
complex projective space $(\mathbb{CP}^{n},\, \lambda \cdot \sigma_{\textup{FS}})$ for some $\lambda>0$ (depending on the volume of $M$).  

If otherwise $\Theta_{1}=\emptyset$, then $\Theta_{0}=\Delta^{M}$
and there are two cases: when $|\mathcal{U}|=1$ and when $|\mathcal{U}|>1$. 
Suppose first that $|\mathcal{U}|=1$. Then
if $\Sigma$ is the only open simplex in the family $\mathcal{U}_0=\mathcal{U}$
the simplex $\Sigma_{\textup{c}}$
has a face which is not contained in $\partial (\Delta^{M})$,
and therefore using the same argument as earlier in the proof we obtain that 
$\Delta^M \neq \Theta_{0} \cup \Theta_{1}$, which
contradicts the fact that $\mathcal{E}^{\Delta}=\mathcal{U}_{0} \cup \mathcal{U}_{1}$
is a coherent family. So $|\mathcal{U}|=1$ may not happen.
Therefore we can pick two 
different simplices $\Sigma_i \in \mathcal{U}_{0}$, $i=0,\,1$.
The statement of the lemma follows.
\end{proof}

From this point on, throughout this proof, we assume that the family $\mathcal{E}^{\Delta}$
contains at least two simplices, since the case where $\mathcal{E}^{\Delta}$ consists
of precisely one simplex is solved in Lemma \ref{firstlemma}.

In what follows let $\widehat F_{\Sigma_i}$, $i=0,\,1$, be the only
$(n-1)$\--dimensional face of $(\Sigma_i)_{\textup{c}}$ which is not
contained in $\partial (\Delta^M)$.

\vspace{2mm}
{\bf Step 2}. 
We give a formula for $\widehat F_{\Sigma_0}$ as a disjoint union
of two subpolytopes of $\Delta^M$, one of which intersects  $\widehat{F}_{\Sigma_0}$
at dimension $n-1$ while the other intersects it at dimension $< n-1$.

\begin{lemma} \label{ll}
For each $x \in \widehat{F}_{\Sigma_0}$
there exists $\Sigma' \in \mathcal{E}^{\Delta}$
such that $x \in (\Sigma')_{\textup{c}}$ and $\Sigma' \neq \Sigma_0$.
\end{lemma}

\begin{proof}
First notice that there exists a unique hyperplane 
$H_{\Sigma_i}$ in $\mathbb{R}^{n}$ which contains $\widehat F_{\Sigma_i}$,
$i=0,\,1$. For each positive integer $n$ let 
$$
U_n:=\mathbb{B}_{1/n}(x) \cap (\Delta^M \setminus (\Sigma_0)_{\textup{c}}),
$$
where $\langle \cdot, \, \cdot \rangle$ denotes the standard scalar
product in $\mathbb{R}^n$, and for each
$x \in \mathbb{R}^n$, $\delta>0$ we write
$\mathbb{B}_{\delta}(x)$ for the standard
open ball in $\mathbb{R}^n$ centered at $x$ and of radius $\delta$ with respect
to $\langle \cdot,\, \cdot \rangle$. We claim that $U_n \neq \emptyset$ for all $n$. Indeed, write
$
H_{\Sigma_0}=\{ x \colon \langle x,\, v \rangle =\lambda \}
$
for some vector $v$ in $\mathbb{R}^n$ and some constant $\lambda \in \mathbb{R}$, 
and suppose that $\Delta^M \subset H_{\Sigma_0}^{+}$
or $\Delta^M \subset H_{\Sigma_0}^{+}$, where $H_{\Sigma_0}^{\pm}$
denote the closed subspaces of $\mathbb{R}^n$ at both sides
of $H_{\Sigma_0}$. Recall the following fact: 

\vspace{1mm}
\emph{Generic fact.} Two subsets $C_1$ and $C_2$ of $\mathbb{R}^n$ are \emph{separated
by a hyperplane $H$} if each lies in a different closed half--space $H^{\pm}$.
If $y$ belongs to the closure of $C_1$, a hyperplane that separates
$C_1$ and $\{y\}$ is called a \emph{supporting hyperplane of $C_1$ at $y$}.
In this case it is a generic (and easy to see) fact that $H \cap \textup{Int}(C_1)=\emptyset$.
\vspace{1mm}

Then $H_{\Sigma_0}$ is a supporting hyperplane
of $\Delta^M$ with respect to any point which is in $\widehat F_{\Sigma_0}$
for $\Delta^M$ and therefore $H_{\Sigma_0} \cap \textup{Int}(\Delta^M)=\emptyset$,
which contradicts $\textup{Int}_{\widehat F_{\Sigma_0}}(\widehat F_{\Sigma_0}) \subset \textup{Int}(\Delta^M)$
(the fact that $\textup{Int}_{\widehat F_{\Sigma_0}}(\widehat F_{\Sigma_0}) \subset \textup{Int}(\Delta^M)$
follows from convexity). Therefore
$\Delta^M \not\subset H_{\Sigma_0}^+$ and $\Delta^M \not\subset H_{\Sigma_0}^-$, which
then implies the existence of $z_i \in \textup{Int}(H^{\pm}_{\Sigma_0})$.
Since $\Delta^M$ and $\mathbb{B}_{1/n}(x)$ are convex, their intersection $\Delta^M \cap \mathbb{B}_{1/n}(x)$
is convex and so we may pick $\epsilon>0$ small enough such that
$$
y_i:=(1-\epsilon)\, x+\epsilon \, z_i \in \Delta^M \cap \mathbb{B}_{1/n}(x), \, \, i=0, \,1.
$$
Since $\langle v, \, z_0 \rangle > \lambda$ and $\langle v, \, z_1 \rangle < \lambda$,
a computation then gives $\langle v, \, y_0 \rangle > \lambda$ and $\langle v, \, y_1 \rangle < \lambda$,
so precisely one of $y_0, \, y_1$ lies in $(\Sigma_0)_{\textup{c}}$ while the other lies in $U_n$, so $U_n \neq \emptyset$
as we wanted to show. For each integer $n$, pick $y_n \in U_n$, and observe that by construction the sequence $\{y_n\}_{n=1}^{\infty}$ 
converges to $x$. Since $\mathcal{E}^{\Delta}$ is finite, there exists a convergent subsequence $\{y_{n_k}\}_{k=1}^{\infty}$ 
of $\{y_n\}_{n=1}^{\infty}$, and a simplex $\Sigma' \in \mathcal{E}^{\Delta}$ such that $y_{n_k} \in (\Sigma')_{\textup{c}}$ for all $k \ge 1$.
Now $\Sigma'\neq \Sigma_0$ because $y_{n_k} \notin \Sigma_0$ but $y_{n_k} \in \Sigma'$, $k \ge 1$,
by construction. Finally since $(\Sigma')_{\textup{c}}$ is compact, $x \in (\Sigma')_{\textup{c}}$ as we wanted
to show. 
\end{proof}

\begin{cor}
Let $\mathcal{F}^{\Delta}$ be the sufamily of $\mathcal{E}^{\Delta}$
consisting of those simplices $\Sigma$ such that
both $\widehat F_{\Sigma_0} \cap \widehat F_{\Sigma} \neq \emptyset$
and $\textup{dim}(\widehat F_{\Sigma_0} \cap \widehat F_{\Sigma})<n-1$,
and $(\mathcal{F}^{\Delta})'$ the subfamily of $\mathcal{E}^{\Delta}$
consisting of those simplices $\widehat{\Sigma} \neq \Sigma_0 $ such that
$\textup{dim}(\widehat F_{\Sigma_0} \cap \widehat F_{\widehat{\Sigma}})=n-1$.
Then $\widehat F_{\Sigma_0}$ may be expressed as the following
union of two subsets of $\Delta^M$:
\begin{eqnarray}\label{mg}
\widehat F_{\Sigma_0}=(\widehat F_{\Sigma_0} \cap (\bigcup_{\Sigma \in \mathcal{F}^{\Delta}}, \Sigma_{\textup{c}}))
\cup (\widehat F_{\Sigma_0} \cap (\bigcup_{\widehat{\Sigma} \in (\mathcal{F}^{\Delta})'} \widehat{\Sigma}_{\textup{c}}))
\end{eqnarray}
and the union is a disjoint one.
\end{cor}

\begin{proof}
The union given in expression (\ref{mg}) is clearly disjoint and we
only need to show that
\begin{eqnarray}\label{mgh}
\widehat F_{\Sigma_0} \subset (\widehat F_{\Sigma_0} \cap (\bigcup_{\Sigma \in \mathcal{F}^{\Delta}} \Sigma_{\textup{c}}))
\cup (\widehat F_{\Sigma_0} \cap (\bigcup_{\widehat{\Sigma} \in (\mathcal{F}^{\Delta})'} \widehat{\Sigma}_{\textup{c}})),
\end{eqnarray}
since the reverse inclusion is trivially true. Notice that showing that 
expression (\ref{mgh}) holds is equivalent to showing that
$$
\widehat F_{\Sigma_0} \subset (\bigcup_{\Sigma \in \mathcal{F}^{\Delta}} \Sigma_{\textup{c}})
\cup (\bigcup_{\widehat{\Sigma} \in (\mathcal{F}^{\Delta})'} \widehat{\Sigma}_{\textup{c}}),
$$
expression which is precisely equivalent to the statement of Lemma \ref{ll},
which concludes the proof.
\end{proof}

\vspace{2mm}

{\bf Step 3}.
We prove that formula \textup{(\ref{mg})}
implies \textup{(}by coherence of $\mathcal{E}^{\Delta}$\textup{)}
that for all $\Sigma \in \mathcal{E}^{\Delta}$, the only faces
of $\Sigma_0$ and of $\Sigma$ which are not contained in $\Delta^M$
are identical, i.e. we have the following.

\begin{lemma} \label{3.5}
$\widehat F_{\Sigma_0}=\widehat F_{\Sigma}$ for all $\Sigma \in \mathcal{E}^{\Delta}$.
\end{lemma}

\begin{proof}
Since $\mathcal{E}^{\Delta}$ is a coherent family
of pairwise disjoint open simplices,
every $(n-1)$\--dimensional
face of every closed simplex $\Sigma_{\textup{c}} \in \mathcal{E}_{\textup{c}}^{\Delta}$, but the
face $\widehat F_{\Sigma}$, is contained in $\partial (\Delta^M)$,
and $\textup{Int}_{\widehat F_{\Sigma}}(\widehat F_{\Sigma}) \subset 
\text{Int}(\Delta^M)$. Therefore
$\widehat F_{\Sigma_0} \cap \Sigma_{\textup{c}}=\widehat F_{\Sigma_0} \cap \widehat F_{\Sigma}$
for all $\Sigma \in \mathcal{E}^{\Delta}$, which by expression (\ref{mg}) then 
implies that:
\begin{eqnarray}\label{mg2}
\widehat F_{\Sigma_0}=\bigcup_{\Sigma \in \mathcal{F}^{\Delta}} (\widehat F_{\Sigma_0} \cap \widehat F_{\Sigma})
\cup
\bigcup_{\widehat{\Sigma} \in (\mathcal{F}^{\Delta})'} (\widehat F_{\Sigma_0} \cap \widehat F_{\widehat{\Sigma}}).
\end{eqnarray}
Let us assume by contradiction that $(\mathcal{F}^{\Delta})'=\emptyset$ and notice that
the left--most member 
of the right hand side of expression (\ref{mg2})
is a union of convex polytopes of dimension
strictly less that $n-1$; furthermore
since $\mathcal{F}^{\Delta}$ is a subfamily of the coherent
family $\mathcal{E}^{\Delta}$, by Remark \ref{2.9} $\mathcal{F}^{\Delta}$ is finite,
and therefore we conclude that $\widehat F_{\Sigma_0}$
is a finite union of convex polytopes, the dimension
of each of which is, by construction of $\mathcal{F}^{\Delta}$,
strictly less than $n-1$, which
is a contradiction since by definition of $\Sigma_0$
we have that $\textup{dim}(\widehat F_{\Sigma_0})=n-1$; 
here we are using the following generic property
of polytopes in $\mathbb{R}^n$: 

\vspace{1mm}
\emph{Generic fact.} \emph{Let $\Delta_0,\, \Delta_1,\ldots, \Delta_k$
be a finite family of polytopes in $\mathbb{R}^n$
such that $\Delta_0=\bigcup_{i=1}^k \Delta_i$
Then there exists $j$ with $1 \le j \le k$ such that
$\textup{dim}(\Delta_0)=\textup{dim}(\Delta_j)$.}
\vspace{1mm}

Therefore
$(\mathcal{F}^{\Delta})'\neq \emptyset$ and hence
there exists $\Sigma^1 \in \mathcal{E}^{\Delta}$
such that both $\Sigma^1 \neq \Sigma_0$ and $\textup{dim}(\widehat F_{\Sigma_0} \cap \widehat F_{\Sigma^1})=n-1$.
Without loss of generality we may assume that $\Sigma_1=\Sigma^1$.
By definition of $H_{\Sigma_i}$,
$
\widehat F_{\Sigma_0} \cap \widehat F_{\Sigma_1} \subset H_{\Sigma_0} \cap H_{\Sigma_1}
$
which in particular implies that 
$$
n-1=\textup{dim}(\widehat F_{\Sigma_0} \cap \widehat F_{\Sigma_1}) \le \textup{dim}(H_{\Sigma_0} \cap H_{\Sigma_1}) \le 
\textup{dim}(H_{\Sigma_0})=n-1$$
and therefore must have $H_{\Sigma_0}=H_{\Sigma_1}$ -- here we are using:

\vspace{1mm}
\emph{Generic fact.} If $L,L'$ are two hyperplanes in $\mathbb{R}^n$ whose intersection
is $(n-1)$--dimensional, then $L=L'$.
\vspace{1mm}

Since $\Delta^M \cap H_{\Sigma_i} =\widehat F_{\Sigma_i}$ and $H_{\Sigma_0}=H_{\Sigma_1}$ we must have
$
\widehat F_{\Sigma_0}=\widehat F_{\Sigma_1}
$.
\end{proof}

\vspace{2mm}

{\bf Step 4}. 
Recall that from Step 3
onwards we have been assuming that the the coherent family $\mathcal{E}^{\Delta}$ contains at least two simplices $\Sigma_0, \, \Sigma_1$. Next we show that the fact that $\Delta^M$ is a Delzant polytope implies that $\Delta^M$ equals the union of $\Sigma_0, \, \Sigma_1$, and hence there are no other simplices in the coherent family $\mathcal{E}^{\Delta}$.

\begin{lemma} \label{lastlemma}
$\mathcal{E}^{\Delta}$
contains precisely two simplices $\Sigma_0$, $\Sigma_1$ joined
at their unique face $(\Sigma_{0})_{\textup{c}} \cap (\Sigma_{1})_{\textup{c}}$, and
$\Delta^{M}=(\Sigma_0)_{\textup{c}} \cup (\Sigma_1)_{\textup{c}}$, where $(\Sigma_0)_{\textup{c}}$
and $(\Sigma_1)_{\textup{c}}$ are joined at their unique common face $\widehat F_{\Sigma_0}$.
\end{lemma}

\begin{proof}
Let us assume that there exists $\Sigma' \neq \Sigma_0,\, \Sigma_1$
with $\Sigma' \in \mathcal{E}^{\Delta}$. By assumption
$\mathcal{E}^{\Delta}$ is a family of pairwise disjoint open simplices,
so by Lemma \ref{3.5} we have that $(\Sigma_0)_{\textup{c}} \cap (\Sigma_1)_{\textup{c}}=(\Sigma_0)_{\textup{c}} \cap (\Sigma')_{\textup{c}}
=(\Sigma_1)_{\textup{c}} \cap (\Sigma')_{\textup{c}}=\widehat{F}_{\Sigma_0}$.
On the other hand, since $\widehat{F}_{\Sigma_0} \subset \partial (\Sigma_0) \cap \partial (\Sigma_1) \cap \partial (\Sigma')$
by construction, we have that 
$$
\emptyset=\textup{Int}((\Sigma')_{\textup{c}}) \cap \textup{Int}((\Sigma_0)_{\textup{c}} \cup (\Sigma_1)_{\textup{c}}), \label{tt2}
$$
and therefore taking the closure of both sides of this expression
we obtain
$$
\emptyset=(\Sigma')_{\textup{c}} \cap 
((\Sigma_0)_{\textup{c}} \cup 
(\Sigma_1)_{\textup{c}})=((\Sigma')_{\textup{c}} \cap (\Sigma_0)_{\textup{c}}) \cup ((\Sigma')_{\textup{c}} \cap 
(\Sigma_1)_{\textup{c}}).
$$
Therefore $\widehat F_{\Sigma_0}=\emptyset$, which is a contradiction. Hence
there does not exist such $\Sigma'$, which then implies that 
$\mathcal{E}^{\Delta} =\{\Sigma_0,\Sigma_1\}$, which proves the first claim of the lemma.

Therefore, since the family $\mathcal{E}^{\Delta}$ is coherent, and by assumption it gives a perfect and equivariant
symplectic packing of $M$, it follows that $\Delta^{M}=(\Sigma_0)_{\textup{c}} \cup (\Sigma_1)_{\textup{c}}$, where $(\Sigma_0)_{\textup{c}}$
and $(\Sigma_1)_{\textup{c}}$ are joined at their unique common face $\widehat F_{\Sigma_0}$ which is in the 
interior of $\Delta^{M}$. 
\end{proof}

\vspace{2mm}

{\bf Step 5}. In this step we analyze which Delzant polytopes in 
$\mathbb{R}^{n}$ may be obtained as the union of the closures of two open simplices.

\begin{figure}
\begin{center}
\epsfig{file=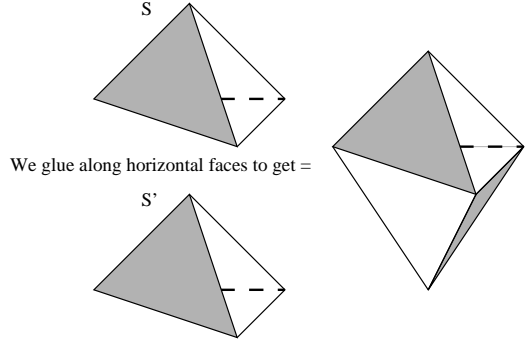}
\label{Picture5}
\caption{A $3$\--dimensional polytope obtained
by gluing two $3$\--dimensional simplices $\textup{S}$ and $\textup{S}\textup{'}$ along a face
does not satisfy the Delzant condition at those vertices contained in the hyperplane 
along which they are glued.}
\end{center}
\end{figure}

\begin{lemma} \label{keylemma}
Let $\Delta$ be a convex polytope in
$\mathbb{R}^{n}$ obtained as the union of the closures of two $n$-dimensional open simplices $\Sigma_{-}$ and $\Sigma_{+}$
joined uniquely by the $(n-1)$\--dimensional face $\widehat F=(\Sigma_{-})_{\textup{c}} \cap (\Sigma_{-})_{\textup{c}}$, 
whose relative interior is contained in the interior of $\Delta$. Then $\Delta$ is unique, and if $n>2$,
$\Delta$ is not a Delzant polytope, i.e. there does not exist a $2n$--dimensional Delzant manifold $M$ 
such that $\Delta=\Delta^{M}$. If $n=1,2$, $\Delta$ is a Delzant polytope if and only
if $\Sigma_{-}$ and $\Sigma_{+}$ are integral simplices.
\end{lemma}
\begin{proof}
First observe that $\Delta$ is unique because $\Delta \subset \mathbb{R}^n$, and $\Delta$
is $n$\--dimensional, and therefore the face $\widehat{F}$ is $(n-1)$\--dimensional, and the plane in which $\widehat{F}$ is contained is uniquely determined by any of its orthogonal vectors, see Remark
\ref{newremark}. 

By assumption $\Delta$ is equal to the union $(\Sigma_{-})_{\textup{c}} \cup (\Sigma_+)_{\textup{c}}$
of the closures of the simplices $\Sigma_{-}$ and $\Sigma_{+}$,
the intersection of which equals an $(n-1)$\--dimensional simplex $\widehat{F}: = (\Sigma_{-})_{\textup{c}} \cap (\Sigma_+)_{\textup{c}}$,
and hence $\widehat{F}$ has precisely $n$ vertices.
Every vertex of $(\Sigma_-)_{\textup{c}}$ and $(\Sigma_+)_{\textup{c}}$ 
outside of $\widehat F$ is also a vertex of $\Delta$. Each other vertex 
of $\widehat F$ is also a vertex of $\Delta$, but the converse need 
not be the case -- if the converse holds, then $\Delta$ has precisely $n+2$ 
vertices, see Figure 6, of which precisely $n$ vertices belong to $\widehat{F}$, and of the two remaining
vertices, one belongs to the simplex $\Sigma_{-}$ but does not belong to $\widehat{F}$, while the other one belongs to the simplex $\Sigma_{+}$ but does not belong to $\widehat{F}$.
Each individual vertex of $\widehat{F}$ is connected with each of the 
other $n+1$ vertices of $\Delta$ by an edge, and this is in 
contradiction with the Delzant property of $\Delta$ unless $n \le 2$, specifically it contradicts the simplicity property that
Delzant polytopes exhibit, c.f. Definition \ref{Definition H}. 

Notice that it follows from Definition \ref{setup} that any interval of finite length which contains precisely one of its two endpoints is an integral simplex. If $n=1$, all three of $\Delta, \, (\Sigma_{-})_{\textup{c}}, \, (\Sigma_{+})_{\textup{c}}$ are closed intervals of finite length and $\Delta$ is the union of $(\Sigma_{-})_{\textup{c}}$ and $(\Sigma_{+})_{\textup{c}}$, the intersection of which consists of 
exactly one point in the interior of $\Delta$, c.f. Figure 3. If $n=2$, see Figure 2, where the packings are explicitly presented. The integrality of $\Sigma_{-}$ and $\Sigma_{+}$ is an essential requirement in order for $\Delta$ to satisfy Definition \ref{Definition H}.
\end{proof}

\begin{remark} \label{newremark}
\normalfont
If in the statement of Lemma \ref{keylemma},
the simplices $\Sigma_{-}$ and $\Sigma_{+}$ where $m$\--dimensional, with $m<n$, there are
infinitely many different ways of gluing them in this fashion; this gluing leads to a convex polytope if and only if $\Sigma_{-}$ and $\Sigma_{+}$ are contained in the same $m$\--dimensional subspace of 
$\mathbb{R}^n$.
\end{remark}

\vspace{2mm} 

{\bf Step 6}. This is the conclusion step. A combination of the previous lemmas gives the proof of Theorem \ref{mt} and Proposition \ref{tt}. Write $X_{n,\,\lambda}=(\mathbb{CP}^n,\, \lambda \cdot \sigma_{\textup{FS}})$ and $Y_{\lambda}=(\mathbb{CP}^1 \times \mathbb{CP}^1, \, \lambda \cdot \sigma_{\textup{FS}})$.

\vspace{2mm}

\emph{Proof of Proposition \ref{tt}}.
Clearly $X_{1, \, \lambda}$ may be packed either by one $2$\--ball or by two $2$-balls
by prescribing a point in its momentum polytope, which is
an interval, c.f. Example \ref{Example B}, whose length depends on the real parameter $\lambda$, and it cannot be packed in any other way, c.f. Figure 3. Similarly,
$Y_{\lambda}$ may be perfectly packed by two equivariant symplectic $2$\--balls, c.f. Figure 2.

Now we show that $X_{n, \, \lambda}$,
$n>2$, may not be packed by two or more balls; if otherwise, there exists a coherent
family of at least two balls, which realizes the perfect packing, and let us call $\Sigma_0$ and $\Sigma_1$ to the corresponding simplices to these
two $n$\--balls. Let $\Delta$ be the momentum polytope of $X_{n,\,\lambda}$ (a simplex in $\mathbb{R}^{n}$) . By Lemma \ref{lastlemma},
$\Delta=(\Sigma_0)_{\textup{c}} \cup (\Sigma_1)_{\textup{c}}$, and hence by Lemma \ref{keylemma},  $\Delta$ is not a Delzant polytope because we are assuming that $n>2$, which is a contradiction. 

It is left to show that $Y_{\lambda}$
may only be perfectly packed by two equivariant symplectic balls, and this is in precisely two ways. Recall that the existence part
is clear. The fact that there are precisely two ways to use two balls to perfectly pack $Y_{\lambda}$
is also clear from Figure 2, since the ball images are integral simplices, see Definition \ref{setup}. Notice that
one equivariantly symplectically embedded ball fills up at most half of the volume of $Y_{\lambda}$,
and therefore it does not give a perfect packing, c.f. Figure 2. Now by Lemma \ref{lastlemma}, $Y_{\lambda}$
does not admit a perfect packing by three or more balls, which concludes the proof. 

\vspace{2mm}

\emph{Conclusion of the proof of Theorem \ref{mt}}.
Lemma \ref{lastlemma}, Lemma \ref{keylemma} and Theorem \ref{Theorem I}, imply that the only 
symplectic--toric manifolds which admit a perfect packing are $X_{n,\,\lambda}$ for arbitrary $n  \ge 1$ and 
$Y_{\lambda}$ for $\lambda>0$. 

The sufficiency condition is implied
by Proposition \ref{tt}, which concludes the proof of Theorem \ref{mt}. 

$\Box$

\section{A Remark on blowing up}

The connection between symplectic ball embeddings and blowing up was first explored by D. McDuff in \cite{M1}. 
Let us first outline McDuff's construction and afterwards we will state the blow up version of
Theorem 1.1 in \cite{P}. Let $(M,\, \sigma)$ be a Delzant manifold and let $J$ be a
$\sigma$--tamed almost complex structure on $M$. 
Recall that we say
that $\sigma$ is $J$\--standard near $p \in M$
if the pair $(\sigma,\, J)$ is diffeomorphic to the standard pair
$(\sigma_0,\, \sqrt{-1})$ of $\mathbb{R}^{2n}$ near $0$.
Choose a $\sigma$\--standard almost complex structure $J$ for which
$\sigma$ is $J$\--standard near $p$ and denote by
$
\Theta:(\widetilde{M},\, \widetilde{J}) \to (M,\, J)
$
to the complex blow up of $M$ at $p$. 
Let $f$ be a
symplectic embedding from $\mathbb{B}_r$ into $M$ which is
holomorphic with respect to the standard multiplication by $\sqrt{-1}$ on
$\mathbb{B}_r$ and $J$, near $0$, and such that
$f(0)=p$. Such an embedding gives rise to a
symplectic form $\widetilde{\sigma}_{f}$ on $\widetilde{M}$ which lies
in the cohomology class
$
[\Theta^* \sigma]-\pi r^2e
$
where $e$ is the Poincar\'e dual to the homology class of the
exceptional divisor $\Theta^{-1}(p)$. The form
$\widetilde{\sigma}_{f}$ is called the {\it symplectic blow up of
$\sigma$ with respect to $f$}, and is uniquely determined up to
isotopy of forms. For the specific construction see \cite{MS} pages
223-225. 

McDuff and Polterovich showed that the same construction extends
for arbitrary symplectic embeddings from $\mathbb{B}_r$ into $M$
without having to assume holomorphicity near $0$. Roughly
speaking, one perturbs the embedding slightly to make it
holomorphic near $0$, and define its blow up as the blow up of 
the perturbed embedding. They also showed that
the isotopy class of the form $\widetilde{\sigma}_{f}$
depends only on the embedding $f$ and the germ of $J$ at
$p$ and that if two symplectic embeddings $f_1$ and $f_2$
are isotopic through a family of symplectic embeddings of
$\mathbb{B}_r$ which take $0$ to $p$, then the
corresponding blow up forms are isotopic. 
On the other hand a (more general) version of the following result was proved in [21]:

\begin{theorem} \label{Theorem L} \textup{(\cite{P})}.
Let $M$ be a $2n$--dimensional Delzant manifold. Then
the space of equivariant symplectic embeddings from 
the $2n$--ball $\mathbb{B}_r$ into $M$ which send the 
origin to the same fixed point $p \in M$,
is homotopically equivalent to the
$n$--torus $\mathbb{T}^n$.
\end{theorem}

And from this result we are able to describe equivariant blow up at a fixed point.

\begin{cor} \label{st1} Let
$f_1$ and $f_2$ be equivariant symplectic
embeddings from the $2n$--dimensional ball $\mathbb{B}_r$ into a $2n$--dimensional Delzant manifold $M$. 
If the normalization condition $f_1(0)=f_2(0)=p$ holds, then
the corresponding blow up manifolds
$(\widetilde{M},\widetilde{\sigma}_{f_1})$ and
$(\widetilde{M},\widetilde{\sigma}_{f_2})$ at $p$ are isotopic, in the sense that the symplectic forms
$\widetilde \sigma_{f_{1}}$ and $\widetilde \sigma_{f_{2}}$ may be joined by a continuous path
$\sigma_{t}$, $0 \le t \le 1$, of symplectic forms with $\sigma_{0}= \widetilde \sigma_{f_{1}}$
and  $\sigma_{1}= \widetilde \sigma_{f_{2}}$.
\end{cor}

\begin{proof}
Follows by observing that if $f_1$ and $f_2$ are equivariant symplectic
embeddings from the $2n$--dimensional ball $\mathbb{B}_r$ into a $2n$--dimensional Delzant manifold $M$
such that $f_1(0)=f_2(0)=p$, then by Theorem \ref{Theorem L}
they are isotopic.
\end{proof}

\section{Further questions}

The following questions regard generalizations of the work presented
in this paper.

\begin{question}
\normalfont
Given $0 \le r \le 1$, find
all Delzant manifolds $M$ such that $\Omega(M)=r$. For $r=0$ it is a trivial
question, and we answered the case $r=1$ in Theorem \ref{mt}. 
What can we say for $r=1/2$? In other words,
to what extent does $\Omega$ encodes the geometry of a Delzant manifold? 
Are there special values of $r$ other
than $0$ and $1$ for which the list of Delzant manifolds $M$ such that
$\Omega(M)=r$ is finite, up to equivariant symplectomorphism?
\end{question}

In \cite{P} we discussed on the topology of the space of partially
equivariant embeddings and suggested a result in this direction.
Recall from \cite{P} that the notion of $\Lambda$--equivariance ($\Lambda \in \mbox{Aut}(\mathbb{T}^n)$)
has a natural extension:

\begin{definition} \label {Definition L}
An embedding
from $\mathbb{B}_r$ into $M$ is \emph{$\sigma$\--equivariant} if there
is a monomorphism $\gamma:\mathbb{T}^{n-k} \to \mathbb{T}^n$,
$k \in \{1,\ldots,n-1 \}$ such that the following diagram commutes
\begin{eqnarray} 
\xymatrix{ \ar @{} [dr] |{\circlearrowleft}
\mathbb{T}^{n-k} \times \mathbb{B}_r  \ar[r]^{\gamma \times f}      \ar[d]^{ \mathbf{\cdot} }  &  \mathbb{T}^{n} \times M 
                  \ar[d]^{\psi}   \\
                   \mathbb{B}_r  \ar[r]^f   &       M    } \nonumber
\end{eqnarray} 
\end{definition}

$M^{\gamma}$ is the set of $p \in M$ such that $\psi(\gamma(t),p)=p$
for all $t \in \mathbb{T}^{n-k}$, and the rest of terminology
is also analogous to that of \cite{P}. 
 
\begin{definition}
We define the real--valued mapping $\Omega_{n-k}$ on the space of $2n$--dimensional
Delzant manifolds by
$
\Omega_{n-k}(M):=(\textup{vol}_{\sigma} (M))^{-1}\, \,  \sup\sum_{\mathcal{E}_{n-k} \in \mathfrak{F}_{n-k}}\textup{vol}_{\sigma} (f(\mathbb{B}_{r}))$
where $\mathfrak{F}_{n-k}$ is the family of sets $\mathcal{E}_{n-k}$ of partially equivariant embeddings of degree $n-k$
and such that their images $f(\mathbb{B}_{r})$ are pairwise disjoint (degree $n-k$ in the sense of the 
monomorphism $\gamma$ defined at the begining of this section having domain a $(n-k)$\--torus)
We say that $M$ admits a \emph{perfect equivariant and symplectic ball packing of degree $n-k$} if $\Omega_{n-k}(M)=1$
at certain family $\mathcal{E}^0_{n-k}$.
\end{definition}

\begin{question} \label{c}
\normalfont
Does for every natural number $n$, every $2n$--dimensional Delzant manifold admit perfect ball packing of degree $k$ for 
all $k$ such that $1 <  k \le n$? Or suppose that $M$ admits a perfect symplectic ball packing, is there $k>1$
such that $M$ admits a perfect equivariant symplectic ball packing of degree $n-k$?
\end{question}

Question \ref{c} is directly related to the partially equivariant version
of Theorem 1.2, which was introduced in \cite{P}.

Finally, Delzant's theorem has been recently generalized in \cite{DuPe} to ``symplectic manifolds whose
principal torus orbits are coisotropic'', and one could ask the 
same question treated in the present paper for the symplectic manifolds in \cite{DuPe}.

\vspace{3mm}
\emph{Acknowledgments.} 
The author is grateful to J.J. Duistermaat  for discussions, and for his hospitality while the author visited Utrecht University  (December 2004 and May 2005). He is also thankful to  D. Auroux and V. Guillemin for discussions and for hosting him at the MIT in a number of occasions (during Fall 2003, Winter 2004).  Finally, he is grateful E. Veomett for comments on a preliminary version, and the referee of the paper for 
his/her suggestions, which have improved the paper.

\noindent
A. Pelayo\\
Department of Mathematics, University of Michigan\\
2074 East Hall, 530 Church Street, Ann Arbor, MI 48109--1043, USA\\
e\--mail: apelayo@umich.edu

\end{document}